\documentclass[11pt]{article}

\usepackage{tikz-cd,amssymb,amsmath,amsthm,amsfonts,sectsty,graphicx,MnSymbol,algorithm,listings,color,appendix,mathrsfs}
\usepackage[letterpaper,hmargin=1.0in,vmargin=1.0in]{geometry}
\newtheorem{theorem}{Theorem}[section]
\newtheorem{lemma}[theorem]{Lemma}

\theoremstyle{definition}

\theoremstyle{remark}

\makeatletter
\def\BState{\State\hskip-\ALG@thistlm}
\makeatother

\graphicspath{ {C:/Users/Bharat/Desktop/LatexImages/} }

\def\Fq{\ensuremath{\mathbb{F}_{\small{q}}}}

\def\R{\ensuremath{\mathbb{R}}}

\def\C{\ensuremath{\mathbb{C}}}

\def\N{\ensuremath{\mathbb{N}}}

\def\Z{\ensuremath{\mathbb{Z}}}
\def\Zplus{\ensuremath{\mathbb{Z}_{\small{\geq 0}}}}

\begin{document}
\title{
  \textbf{A combinatorial proof of Bass's determinant formula for the zeta function of regular graphs}}
\author{
  \Large{Bharatram Rangarajan}\footnote{School of Computer Science, Tel Aviv University }
  }
\setlength{\parskip}{1ex plus 0.5ex minus 0.2ex}
\maketitle

\begin{abstract}
We give an elementary combinatorial proof of Bass's determinant formula for the zeta function of a finite regular graph. This is done by expressing the number of non-backtracking cycles of a given length in terms of Chebychev polynomials in the eigenvalues of the adjacency operator of the graph.  
\end{abstract}

\section{Introduction}
In the 158 years since Bernard Riemann published his seminal work "On the Number of Primes Less Than a Given Magnitude" \cite{riemann}, there have been several generalizations of the Riemann zeta function in various settings. Prominent examples include the Dedekind zeta function of algebraic number fields, the Hasse-Weil zeta function of an algebraic variety, and the Selberg zeta function of a hyperbolic surface.\\
Broadly speaking, a zeta function is a complex function which when expressed as an appropriate series, yields a coefficient sequence that counts "objects" of a given "weight" assembled from an underlying set of building blocks or "primes". For instance, the Riemann zeta function corresponds to a Dirichlet series where the coefficient of $1/k^s$ counts the number of positive integers (constructed using the primes of $\Z$ as building blocks) of absolute value $k$ (which in this case is trivially $1$ for every $k \in \N$). Similarly the Hasse-Weil zeta function corresponds to an ordinary power series that counts the number of positive divisors (constructed using the places of the function field acting as "primes").\\
The utility of a zeta function arises from the fact that many interesting properties of the underlying structure can be inferred from the zeros and poles of the corresponding zeta function. For instance, the zeros of the Riemann zeta function are related to the distribution of prime numbers, and the zeros of the Hasse-Weil zeta function of a projective curve over a finite field are related to the number of rational points on the curve. The zeros and poles of the Selberg zeta function appear in the Selberg trace formula, which relates the distribution of primes with the spectrum of the Laplace-Beltrami operator of the surface.\\
The precursor to the zeta function of a graph, as we know it today, is the Selberg zeta function of a Riemannian manifold. For a hyperbolic surface $M=\Gamma \slash H$, the Selberg zeta function $\gamma_{M}(s)$ is an Euler product over the set of all primitive closed geodesics of $M$. The zeros and poles of the Selberg zeta function appear in the Selberg trace formula, which relates the distribution of primes with the spectrum of the Laplace-Beltrami operator of the surface. This line of study was further extended by Ihara \cite{ihara} to obtain a $p$-adic analogue of the Selberg trace formula, opening up further avenues for the study of geodesic zeta functions in discrete settings. 
The idea of considering closed geodesics as primes inspired the work of Hashimoto \cite{hashimoto}, Hyman Bass \cite{iharabass}, Kotani and Sunada \cite{sunada} to come up with an analogous notion in the discrete setting of a finite graph, using the prime cycle classes of the graph in place of primitive geodesics.\\

Formally, for a finite graph $G=(V,E)$, the (Ihara) zeta function of $G$, denoted $\zeta_G(t)$, is defined as the Euler product
$$\prod \limits_{[P] \in \mathcal{P}} \frac{1}{1-t^{|P|}}$$
where $\mathcal{P}$ is the set of primitive cycle classes. These notions are made rigorous in the subsequent section.\\

Just like the Selberg zeta function is related to the spectrum of the Laplace-Beltrami operator of the surface, it is natural to ask if its discrete analogue, the Ihara zeta function of a graph, is related to the spectrum of the Laplacian matrix (or the adjacency matrix) of the graph. This is precisely the result of Bass \cite{iharabass} who gives an elegant expression for the Ihara zeta function of a graph $G=(V,E)$ as the rational function
$$\zeta_G(t) = \frac{1}{(1-t^2)^{|E|-|V|} det(I-tA+(D-I)t^2)}$$
where $A$ is the adjacency matrix of $G$ and $D$ is the diagonal matrix of degrees of the vertices of $G$, or in other words, $D=diag(A\vec{1})$. In particular, if $G$ is $d$-regular, then
$$\zeta_G(t) = \frac{1}{(1-t^2)^{|E|-|V|} det(I-tA+(d-1)t^2I)}$$
which immediately gives us a way of obtaining precisely the set of poles of $\zeta_G(t)$.\\

The significance of the poles of the zeta function arises from a surprising analogue of the classical Riemann hypothesis in our present context. The classical Riemann hypothesis for the Riemann zeta function $\zeta(t)$ states that every non-trivial zero of $\zeta(t)$ lies on the line $Re(z)=1/2$ in the complex plane. Analogues of the Riemann hypothesis can be formulated for other zeta functions too. For instance, the Riemann hypothesis for curves over finite fields states that every zero of the Hasse-Weil zeta function for a projective curve over a finite field $\Fq$ is of absolute value exactly $1/\sqrt{q}$. It is interesting to note that while the classical Riemann hypothesis remains elusive, the Riemann hypothesis for finite fields has been proved, and is one of the crowning achievements of twentieth-century mathematics.\\
It is natural to ask what the appropriate formulation of the Riemann hypothesis is for the Ihara zeta function, and what it means for the graph. A $d$-regular graph $G$ is said to be \emph{Ramanujan} if for every eigenvalue $\mu \in \R$ of the adjacency matrix of $G$ with $|\mu| \neq d$ satisfies
$$|\mu| \leq 2\sqrt{d-1}$$
Combining this with Bass's determinant formula for the zeta function of $G$, if can be easily shown \cite{murty} that
\begin{lemma}
A $d$-regular graph $G$ is Ramanujan iff every pole $\lambda \in \C$ of $\zeta_G(t)$ such that $|\lambda| \neq \pm 1$ and $|\lambda| \neq \pm (d-1)^{-1}$ satisfies
$$|\lambda| = \frac{1}{\sqrt{d-1}}$$
\end{lemma}
Thus the Ramanujan property elegantly reflects in the poles of the Ihara zeta function of the graph, and the expression mirrors the Riemann hypothesis for the Hasse-Weil zeta function of curves over finite fields. For a brief survey of Ramanujan graphs and their significance, the reader is referred to Murty's monograph \cite{murty}.\\

There exist several proofs of Bass's determinant formula \cite{sunada} \cite{starkterras}, and most proofs start by expressing the zeta function in terms of not the adjacency matrix $A$ of $G$, but the adjacency matrix $H$ of the oriented line graph of $G$ (called the Hashimoto edge-incidence matrix). This is followed by appropriate linear-algebraic manipulations of the matrices involved in order to arrive at the desired expression. Foata and Zeilberger \cite{zeil} presented an an insightful combinatorial proof employing the algebra of Lyndon words.

In this paper, we shall see a more elementary combinatorial proof of Bass's determinant formula in the special case when $G$ is regular. While the assumption of regularity is certainly a limitation, it allows for a more transparent and natural proof The basic proof idea is outlined as follows:
\begin{itemize}
\item We observe that the zeta function $\zeta_G(t)$ has an expansion of the form  
$$\zeta_G(t)= exp\left( \sum \limits_{k=1}^{\infty} N_k \frac{t^k}{k} \right)$$
where for $k \in \N$, $N_k$ is the number of rooted, \emph{non-backtracking} cycles in $G$ of length $k$. This is explored in section $2$.
\item While an expression for $N_k$ is not immediate, a natural starting point is the study of non-backtracking walks on $G$. We can construct the family $\{A_k\}_{k \in \Zplus}$ of $n \times n$ matrices such that for every $k \in \Zplus$ and every $v,w \in V$, $(A_k)_{v,w}$ is the number of non-backtracking walks on $G$ of length $k$ from $v$ to $w$. We shall discuss the construction of these non-backtracking walk matrices in section $3$.
\item While it might be tempting to claim that $N_k=Tr(A_k)$, unfortunately that is not the case. However, while they may not be equal, they are indeed precisely related. In section $4$, we develop a combinatorial lemma to relate $N_k$ and $Tr(A_k)$.
\item The combinatorial lemma greatly simplifies the problem since $Tr(A_k)$ is well-understood in terms of the eigenvalues of $A$ and a family of orthogonal polynomials called the Chebychev polynomials. We shall put these ingredients together in section $5$ to conclude with a proof of Bass's determinant formula.
\end{itemize}

\section{Preliminaries}
For an integer $d \geq 2$, let $G=(V,E)$ be a finite $d$-regular undirected graph with adjacency matrix $A$. A \emph{walk} on the graph $G$ is a sequence $v_0v_1\dots v_{k}$ where $v_0,v_1,\dots,v_k$ are (not necessarily distinct) vertices in $V$, and for every $0 \leq i \leq k-1$, $(v_i,v_{i+1}) \in E$. The vertex $v_0$ is referred to as the \emph{root} (or origin) of the above walk, $v_k$ is the terminus of the walk, and the walk is said to have length $k$.\\ 
It is often useful to equivalently define a walk as a sequence of directed or oriented edges. Associate each edge $e=(v,w) \in E$ with two directed edges (or rays) denoted
$$\vec{e}=(v \to w)$$
$$\vec{e}^{-1} = (w \to v)$$
Note that the origin $org(\vec{e})$ is the vertex $v$ and its terminus $ter(\vec{e})$ is the vertex $w$. Similarly, the origin $org(\vec{e}^{-1})$ is the vertex $w$ and its terminus $ter(\vec{e})$ is the vertex $v$. Let $\vec{E}$ denote the set of $m=nd$ directed edges of $G$. So a walk of length $k$ can equivalently be described as a sequence $\vec{e}_1\vec{e}_2\dots\vec{e}_k$ of $k$ (not necessarily distinct) oriented edges in $\vec{E}$ such that for every $1 \leq i \leq k-1$,
$$ter(\vec{e}_i) = org(\vec{e}_{i+1})$$
This is a walk that starts at $org(\vec{e}_1)$ and ends at $ter(\vec{e}_k)$.\\
It is easy to show that for any $k \in \N$, the number of walks of length $k$ between vertices $u,v \in V$ is exactly $(A^k)_{u,v}$. In particular, the total number of rooted cycles of length $k$ in $G$ is exactly
$$Tr(A^k)$$

A \emph{non-backtracking walk} of length $k$ from $v_0 \in V$ to $v_k \in V$ is a walk $v_0v_1\dots v_k$ such that for every $1 \leq i \leq k-1$,
$$v_{i-1} \neq v_{i+1}$$
Equivalently, a non-backtracking walk of length $k$ from $v \in V$ to $w \in V$ is a walk $\vec{e}_1\vec{e}_2\dots\vec{e}_k$
such that $org(\vec{e}_1)=v$, $ter(\vec{e}_k)=w$ and for every $1 \leq i \leq k-1$,
$$\vec{e}_{k+1} \neq \vec{e}^{-1}_{k}$$
Non-backtracking random walks on graphs have been studied in the context of mixing time \cite{alon}, cut-offs \cite{peres}, and exhibit more useful statistical properties than ordinary random walks. In \cite{peres}, the authors obtain further interesting results on the eigendecomposition of the Hashimoto matrix $H$.\\
 
A rooted, non-backtracking cycle of length $k$ with root $v$ is a non-backtracking walk $v,v_1,v_2,\dots,v_{k-1},v$ with the additional boundary constraint that
$$v_1 \neq v_{k-1}$$

Let $\mathcal{C}$ denote the set of all rooted, non-backtracking, closed walks in $G$, and for $C \in \mathcal{C}$, let $|C|$ denote the length of the walk $C$. There are two elementary constructions we can carry out to generate more elements of $\mathcal{C}$ from a given cycle $C$:
\begin{itemize}
\item \emph{Powering}: Given a rooted, non-backtracking closed walk $C \in \mathcal{C}$ of length $k$ of the form
$$C=\vec{e}_1\vec{e}_2\dots\vec{e}_k$$
then for $m \geq 1$ define a power
$$C^m = \underbrace{\vec{e}_1\dots\vec{e}_k \vec{e}_1\dots\vec{e}_k \dots \vec{e}_1\dots\vec{e}_k}_{\text{m times}}$$
which is a concatenation of the string of edges corresponding to the walk $C$ with itself $m$ times. Note that $C^m$ is also a rooted, non-backtracking closed walk in $G$ of length $mk$. Essentially, $C^m$ represents the walk obtained by repeating or winding the walk $C$ $m$ times. Also note that $C$ and $C^m$ are both rooted at the same vertex.
\item \emph{Cycle class}: Given a rooted, non-backtracking closed walk $C \in \mathcal{C}$ of length $k$ of the form
$$C=\vec{e}_1\vec{e}_2\dots\vec{e}_k$$
we can form another walk
$$C^{(2)}=\vec{e}_2 \vec{e}_3 \dots \vec{e}_k \vec{e}_1$$
which is also a rooted, non-backtracking closed walk in $G$ of length $k$, but now rooted at the origin of the directed edge $\vec{e}_2$ (or the terminus of $\vec{e}_1$). More generally, for $1 \leq j \leq k$, define
$$C^{(j)} = \vec{e}_j \vec{e}_{j+1} \dots \vec{e}_k \vec{e}_1 \vec{e}_2 \dots \vec{e}_{j-1}$$
which is a cyclic permutation of the walk $C$ obtained by choosing a different root. So given a walk $C \in \mathcal{C}$ of length $k$, we get $k-1$ additional walks in $\mathcal{C}$ of length $k$ for free this way. In fact, this defines an equivalence class $\sim$ on $\mathcal{C}$, and the set
$$[C]=\{C^{(1)},C^{(2)},\dots,C^{(k)} \}$$
is called the equivalence class of $C$. An element $[C] \in \mathcal{C}/\sim$ represents a non-backtracking closed walk modulo a choice of root.
\end{itemize}
Consider the operation of powering. As mentioned, for an element $C \in \mathcal{C}$ rooted at a vertex $v$, $C^m$ is also rooted at $v$. It is tempting to ask if the notion of powering can be naturally extended to a product of elements of $\mathcal{C}$ as long as we work with a fixed root vertex $v \in V$. For a vertex $v \in V$, let $\mathcal{C}_{v} \subseteq \mathcal{C}$ be the set of non-backtracking closed walks on $G$ rooted at the vertex $v$. Then for two elements $C_1, C_2 \in \mathcal{C}_{v}$, we could try defining the product $C_1C_2$ simply as the concatenation (or composition) of the walks $C_1$ and $C_2$. By this definition, $C_1C_2$ would be a closed walk rooted at $v$, but is not necessarily non-backtracking. For instance, if the last edge of the walk $C_1$ is $\vec{e}$ and the first edge of the walk $C_2$ is $\vec{e}^{-1}$, then the walk $C_1C_2$ clearly has a backtracking instance of the form $\vec{e} \vec{e}^{-1}$.\\
Even though we do not have a simple notion of multiplication of elements of $\mathcal{C}$ (even if we fix a root), we can still try to define a notion of an irreducible or \emph{prime} walk using the available powering operation as follows: a walk $P \in \mathcal{C}$ shall be called a prime walk if there exists no element $C \in \mathcal{C}$ and $m \geq 2$ such that $P=C^m$. Intuitively, a prime walk in $\mathcal{C}$ is one that is not a repeated winding of a simpler closed walk in $\mathcal{C}$. Note that every element of $\mathcal{C}$ is either a prime or a prime power. In particular, the set of primes of $\mathcal{C}$ act as the basic building blocks of $\mathcal{C}$ under the operation of powering.\\ 
Let $\mathcal{P}$ denote the set of equivalence classes of primes. The Euler product
$$\prod \limits_{[P] \in \mathcal{P}} \frac{1}{1-t^{|P|}}$$
is called the \emph{Ihara zeta function} of the graph $G$, denoted $\zeta_G(t)$.\\
Let $N_k$ denote the number of rooted, non-backtracking cycles in $G$ of length $k$. Then observe that
$$\sum \limits_{k=1}^{\infty} N_k \frac{t^k}{k} = \sum  \limits_{\text{prime }P} \frac{1}{|P|} \left( \sum \limits_{m=1}^{\infty} \frac{t^{m|P|}}{m} \right) = - \sum \limits_{[P] \in \mathcal{P}} \log{(1-t^{|P|})}$$
Thus,
$$\zeta_G(t) = \prod \limits_{[P] \in \mathcal{P}} \frac{1}{1-t^{|P|}} = exp \left( \sum \limits_{k=1}^{\infty} N_k \frac{t^k}{k} \right)$$

Just like the number of (rooted) cycles in $G$ of length $k$ is $Tr(A^k)$, we can describe the number $N_k$ of rooted, non-backtracking cycles in $G$ of length $k$ as the trace of the matrix $H^k$ where $H$ is the \emph{Hashimoto edge incidence matrix} of $G$ defined as follows: $H \in \C^{dn \times dn}$ with
$$H_{i,j} = \begin{cases}
1 & \text{ if } \vec{e}_j \neq \vec{e}_i \text{ and }ter(\vec{e}_i) = org(\vec{e}_j)\\
0 & \text{ otherwise}
\end{cases}$$
In other words, the entry $H_{i,j}$ is an indicator for whether the oriented edge $\vec{e}_i$ feeds into the oriented edge $\vec{e}_j$ allowing us to form a non-backtracking walk $\vec{e}_i \vec{e}_j$ of length $2$.\\
It is clear that for every $k \in \mathcal{N}$,
$$N_k = Tr(H^k)$$
and so 
$$\zeta_G(t) = exp \left( \sum \limits_{k=1}^{\infty} Tr(H^k) \frac{t^k}{k} \right) = -Tr\left( \log{(I-tH)} \right)$$
By Jacobi's formula relating the trace of the logarithm of a matrix to the logarithm of its determinant, we get
$$\zeta_G(t) = \frac{1}{det(I-Ht)}$$
In particular, this establishes the rationality of the Ihara zeta function of a regular graph, and further implies that the reciprocal $\zeta_G(t)^{-1}$ is a polynomial in $t$ over $\Z$ of degree at most $m=nd$. However, it is not immediate what the spectrum of $H$ is. Thus, in a sense, Bass's determinant formula can be interpreted as an expression that allows us to determine the spectrum of the Hashimoto matrix $H$ in terms of the spectrum of the adjacency matrix $A$.

\section{Non-backtracking walks and Chebychev polynomials}
Just like $(A^k)_{v,w}$ counts the total number of walks on $G$ from $v$ to $w$ (with backtrackings) of length $k$, we can construct a family
$$A_0,A_1,A_2,A_3,\dots$$
of $n \times n$ matrices over $\C$ such that the value $(A_k)_{v,w}$ is the number of non-backtracking walks on $G$ from $v$ to $w$ of length $k$. This family $\{A_k\}_{k \in \N}$ can be inductively defined using powers of $A$ as follows:
\begin{itemize}
\item $A_0=I$
\item $A_1=A$
\item $A_2=A^2-dI$
\item For $k \geq 3$,
$$A_k=A_{k-1}A-(d-1)A_{k-2}$$ 
\end{itemize}
The recurrence relation above can be used to easily show that the ordinary (matrix) generating function for the above sequence is
$$\sum \limits_{k=0}^{\infty} t^k A_k = (1- t^2)I. \left( I-tA + (d-1)t^2 I \right)^{-1}$$
With some abuse of notation, we can rewrite this generating function as
$$\frac{1-t^2}{1-At+(d-1)t^2}$$

The generating function above is closely related to the generating function of a well-studied family of orthogonal polynomials. Consider the family of Chebychev polynomials of the second kind
$$U_0(x), U_1(x), U_2(x), \dots $$
of univariate complex polynomials defined by the recurrence
$$U_0(x)=1$$
$$U_1(x) = 2x$$
and for $k \geq 2$,
$$U_k(x) = U_{k-1}(x)U_1(x) - U_{k-2}(x)$$
and with generating function
$$\sum \limits_{k=0}^{\infty} U_k(x)t^k = \frac{1}{1-2xt+t^2}$$

It is easy to see that
$$\sum \limits_{0 \leq j \leq k/2} A_{k-2j} = (d-1)^{k/2} U_k \left( \frac{A}{2\sqrt{d-1}} \right)$$
implying that for $k \geq 2$,
$$A_k = (d-1)^{k/2} U_k \left( \frac{A}{2 \sqrt{d-1}} \right) - (d-1)^{k/2-1} U_{k-2} \left( \frac{A}{2 \sqrt{d-1}} \right)$$
In fact this expression can be made to hold consistently for $0 \leq k \leq 2$ too by assigning
$$U_m(x) = 0$$
for every $m < 0$. This allows us to work with the above expression for $A_k$ for \emph{all} non-negative integers $k$.\\
Taking trace on both sides,
$$Tr(A_k) = (d-1)^{k/2} \sum \limits_{j=0}^{n-1} U_k \left( \frac{\mu_i}{2 \sqrt{d-1}} \right) - (d-1)^{k/2-1} \sum \limits_{i=0}^{n-1} U_{k-2} \left( \frac{\mu_i}{2 \sqrt{d-1}} \right) $$
where 
$$d=\mu_0 \geq \mu_1 \geq \dots \geq \mu_{n-1} \geq d$$
are the $n$ eigenvalues of the adjacency matrix $A$. Thus we have an expression for the trace of $A_k$ as a polynomial in the eigenvalues of $A$. This approach is used in the seminal work of Lubotzky, Phillips and Sarnak in their construction of Ramanujan graphs \cite{lps}, and for a detailed and elementary exposition of Chebychev polynomials and non-backtracking walks on regular graphs, the reader is referred to the monograph by Davidoff, Sarnak and Valette \cite{sarnak}. \\

While $(A_k)_{v,w}$ counts the number of walks on $G$ from vertex $v$ to vertex $w$ without backtracking, observe that the diagonal element $(A_k)_{v,v}$ does \emph{not} count the number of non-backtracking cycles of length $k$ rooted at $v$. This is because $(A_k)_{v,v}$ also counts walks of the form
$$e_1 e_2 \dots e_k$$
where $e_{i+1} \neq \overline{e_i}$ for any $1 \leq i \leq k-1$ but $e_k = \overline{e_1}$. That is, $e_1 e_2 \dots e_k$ is non-backtracking as a walk from $v$ to $v$, but when considered as a closed walk (or a loop), the two end edges form a backtracking! Such an instance of a backtracking that gets overlooked in $Tr(A_k)$ shall be referred to as a \emph{tail}. \\
So $Tr(A_k)$ counts the number of closed, rooted walks of length $k$ that could have at most $1$ tail, and hence does \emph{not} count the  rooted, non-backtracking cycles of length $k$. It is interesting to ask what the number of closed, rooted non-backtracking walks of length $k$ is.\\
Denote $Tr(A_k)$ by $M_k$. In the following section, we shall establish a combinatorial lemma relating $M_k$ with $N_k$. 

\section{The Combinatorial Lemma}
Firstly it is clear that
$$N_1=N_2=0$$
just like $M_1=M_2=0$. For $k \geq 3$, observe that we can count the number $M_k$ of \emph{tailed} non-backtracking, closed walks of length $k$ based on the length of the tail as follows: Every tailed non-backtracking closed walks of length $k$ can be constructed as:
\begin{itemize}
\item A tailless, rooted, non-backtracking closed walk of length $k$, and there are $N_k$ of them.\\
\centerline{\includegraphics*[scale=.4]{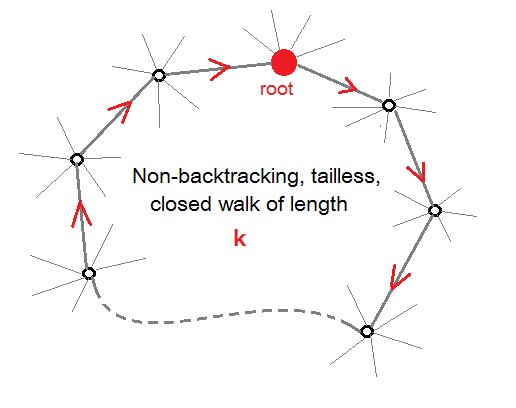}}\\
\item A tailless, rooted, non-backtracking closed walks of length $k-2$ and a tail of length $1$. Since the root is fixed and there are $d-2$ choices for the tail (and consequently, the new root), the number of rooted, non-backtracking closed walks of length $k$ with a tail of length $1$ is $(d-2)N_{k-2}$.\\
\centerline{\includegraphics*[scale=.4]{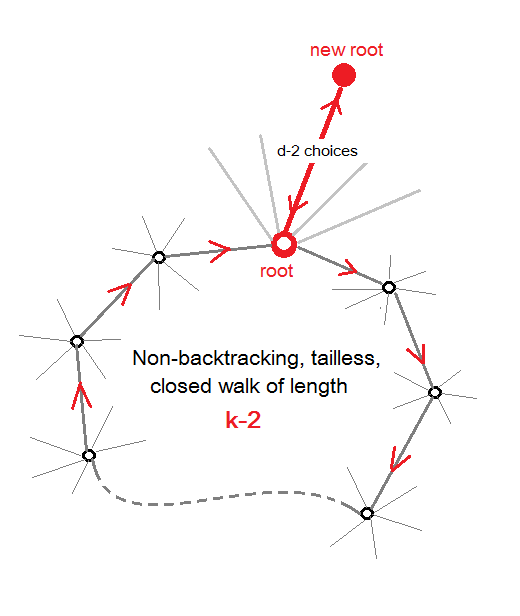}}\\
\item A tailless, rooted, non-backtracking closed walks of length $k-4$ and a tail of length $2$. In this case the first vertex of the tail can be chosen in $d-2$ ways, and the next vertex (the new root) can be chosen in $d-1$ ways. So the number of rooted, non-backtracking closed walks of length $k$ with a tail of length $2$ is $(d-1)(d-2)N_{k-4}$.\\
\centerline{\includegraphics*[scale=.4]{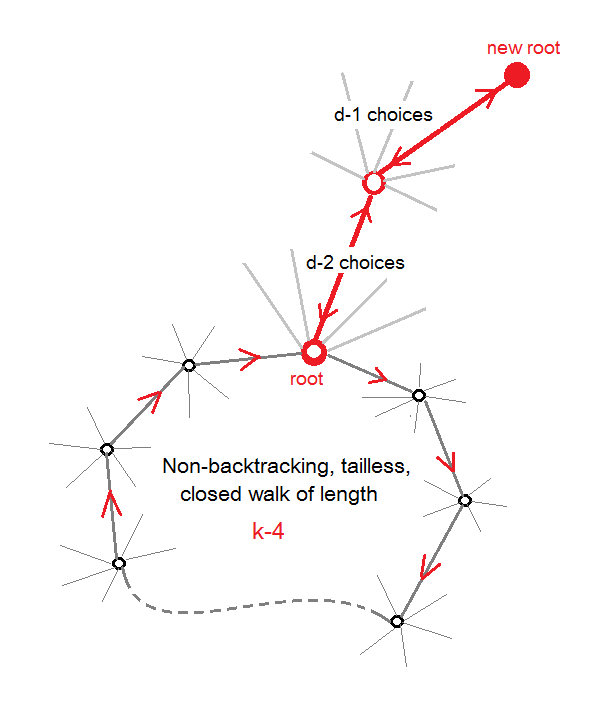}}\\
\item A tailless, rooted, non-backtracking closed walks of length $k-6$ and a tail of length $3$. In this case the first vertex of the tail can be chosen in $d-2$ ways, the next vertex can be chosen in $d-1$ ways, and the new root can be chosen in $d-1$ ways. So the number of rooted, non-backtracking closed walks of length $k$ with a tail of length $3$ is $(d-1)^2(d-2)N_{k-6}$.\\
\centerline{\includegraphics*[scale=.4]{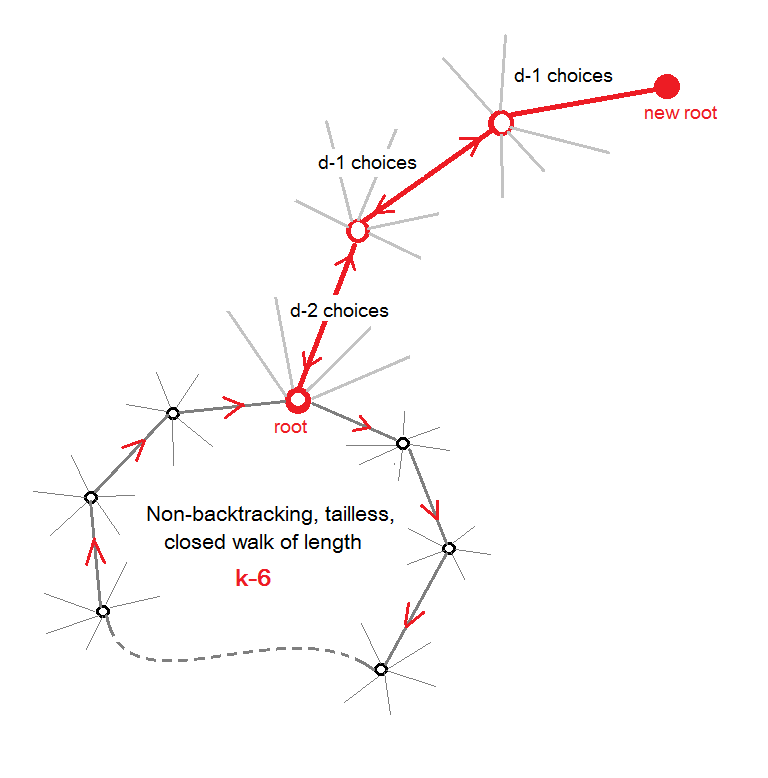}}\\
\item More generally, for $2 \leq r \leq \lfloor k/2 \rfloor$, the number of rooted, non-backtracking closed walks of length $k$ with a tail of length $r$ is $(d-1)^{r-1}(d-2)N_{k-2r}$.\\
\centerline{\includegraphics*[scale=.4]{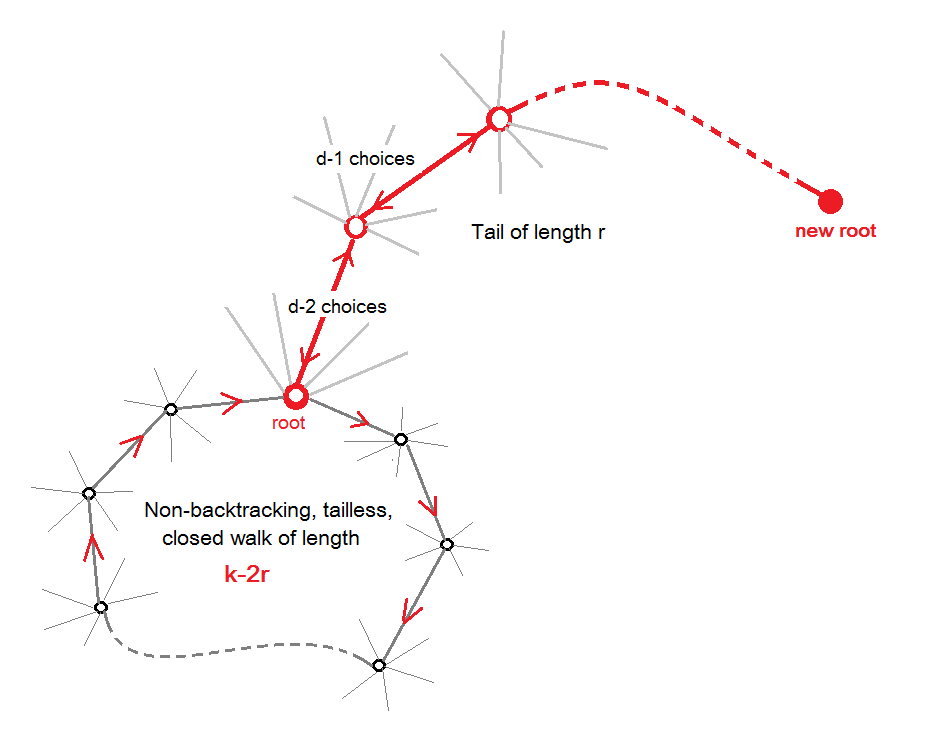}}\\
\end{itemize} 
Thus for every $k \geq 3$,
$$M_k = N_k + (d-2)N_{k-2} + (d-2)(d-1) N_{k-4} + (d-2)(d-1)^2 N_{k-6} + \dots + (d-2)(d-1)^{\lfloor k/2 \rfloor - 1}N_{k-2 \lfloor k/2 \rfloor}$$

While this expression looks cumbersome, a straightforward rearrangement shows that this is equivalent to
\begin{lemma} For every $k \geq 3$,
$$N_k = \begin{cases}
M_k - (d-2) (M_{k-2}+M_{k-4}+\dots+M_1) & \text{ if k is odd}\\
M_k - (d-2) (M_{k-2}+M_{k-4}+\dots+M_2) & \text{ if k is even}\\
\end{cases}$$
\end{lemma}

\section{The Bass determinant formula}
From the combinatorial lemma established in the previous section, we see that
$$N_k=\begin{cases}
Tr(A_k) - (d-2) \left(Tr(A_{k-2})+Tr(A_{k-4}+\dots+Tr(A_1)\right) & \text{ if k is odd}\\
Tr(A_k) - (d-2) \left(Tr(A_{k-2})+Tr(A_{k-4})+\dots+Tr(A_2)\right) & \text{ if k is even}\\
\end{cases}$$
which implies, by linearity of trace, that
$$N_k=\begin{cases}
Tr\left( A_k - (d-2)(A_{k-2}+A_{k-4}+\dots+A_1) \right) & \text{ if k is odd}\\
Tr\left( A_k - (d-2)(A_{k-2}+A_{k-4}+\dots+A_2) \right) & \text{ if k is even}\\
\end{cases}$$
Recall that
$$\sum \limits_{0 \leq j \leq k/2} A_{k-2j}=(d-1)^{k/2} U_k\left( \frac{A}{2\sqrt{d-1}}\right) $$
So for odd $k$
\begin{align*}
A_k - (d-2)(A_{k-2}+A_{k-4}+\dots+A_1) & = (A_k + A_{k-2}+\dots+A_1 ) - (d-1) (A_{k-2}+A_{k-4}+\dots+A_1)\\
& = (d-1)^{k/2} U_k\left( \frac{A}{2\sqrt{d-1}}\right)-(d-1)^{k/2} U_{k-2}\left( \frac{A}{2\sqrt{d-1}}\right)
\end{align*}

Similarly for even $k$,
\begin{align*}
A_k - (d-2)(A_{k-2}+A_{k-4}+\dots+A_2) & = (A_k + A_{k-2}+\dots+A_2 + A_0 ) - (d-1) (A_{k-2}+A_{k-4}+\dots+A_2 + A_0) \\
& = (d-1)^{k/2} U_{k-2}\left( \frac{A}{2\sqrt{d-1}}\right)-(d-1)^{k/2} U_{k-2}\left( \frac{A}{2\sqrt{d-1}}\right) + (d-2)I 
\end{align*}

As it so happens, the polynomial
$$U_k(x) - U_{k-2}(x) = 2 T_k(x)$$
where $T_k(x)$ is called the \emph{Chebychev polynomial of the first kind} of order $k$. The Chebychev polynomials of the first kind are defined in a way very similar to the Chebychev polynomials of the second kind:
$$T_0(x)=1$$
$$T_1(x)=x$$
and for $k \geq 2$,
$$T_k(x) = 2x T_{k-1}(x) - T_{k-2}(x)$$
It is easy to show that $T_k(x)$ has a generating function
$$\sum \limits_{k=0}^{\infty} T_k(x)t^k = \frac{1-xt}{1-2xt+t^2}$$
It is convenient to express $N_k$ in terms of Chebychev polynomials of the first kind as follows:
$$N_k=\begin{cases}
Tr\left( 2(d-1)^{k/2}T_k \left( \frac{A}{2\sqrt{d-1}} \right)  \right) & \text{ if k is odd}\\
Tr\left( 2(d-1)^{k/2} T_k \left( \frac{A}{2\sqrt{d-1}} \right) + (d-2)I \right) & \text{ if k is even}\\
\end{cases}$$
This simplifies to
$$N_k=\begin{cases}
\sum \limits_{j=0}^{n-1} 2(d-1)^{k/2} T_k \left( \frac{\mu_j}{2\sqrt{d-1}} \right) & \text{ if k is odd}\\
n(d-2) + \sum \limits_{j=0}^{n-1} 2(d-1)^{k/2} T_k \left( \frac{\mu_j}{2\sqrt{d-1}} \right)    & \text{ if k is even}\\
\end{cases}$$
The generating function for $N_k$ given by
\begin{align*}
\sum \limits_{k=1}^{\infty} N_k t^k & = n(d-2)(t^2+t^4+t^6+\dots) +  \sum \limits_{k=1}^{\infty} t^k \left( \sum \limits_{j=0}^{n-1} 2(d-1)^{k/2} T_k \left( \frac{\mu_j}{2\sqrt{d-1}} \right) \right) \\
& = n(d-2)(t^2+t^4+t^6+\dots) + 2\sum \limits_{j=0}^{n-1} \sum \limits_{k=1}^{\infty} (t\sqrt{d-1})^k T_k \left( \frac{\mu_j}{2\sqrt{d-1}} \right) \\
& = n(d-2)(t^2+t^4+t^6+\dots) + \sum \limits_{j=0}^{n-1} \left( \frac{2-\mu_j t}{1-\mu_j t + (d-1)t^2}-1 \right)\\
& = n(d-2)\frac{t^2}{1-t^2} + \sum \limits_{j=0}^{n-1} \frac{\mu_j t - 2(d-1)t^2}{1-\mu_j t + (d-1)t^2}
\end{align*}

Thus,
$$N_1 + N_2t + N_3t^2 + \dots = n(d-2) \frac{t}{1-t^2} +  \sum \limits_{j=0}^{n-1} \frac{\mu_j  - 2(d-1)t}{1-\mu_j t + (d-1)t^2} $$
While this expression does not seem very elegant stated this way, observe that the derivative of $1-t^2$ is $-2t$, and the derivative of $1-\mu_j t + (d-1)t^2$ is $-\mu_j + 2(d-1)t$. Rewriting the above expression to highlight this observation, 
$$N_1 + N_2t + N_3t^2 + \dots = -\frac{n(d-2)}{2} \frac{-2t}{1-t^2} -  \sum \limits_{j=0}^{n-1} \frac{-\mu_j  + 2(d-1)t}{1-\mu_j t + (d-1)t^2} $$
This suggests that we could integrate both sides to obtain
\begin{align*}
N_1t + N_2 \frac{t^2}{2} + N_3 \frac{t^3}{3} + \dots & = -\frac{n(d-2)}{2} \log{(1-t^2)} - \sum \limits_{j=0}^{n-1} \log{(1-\mu_j t + (d-1)t^2)} \\ 
& = -\left( \frac{nd}{2} - n \right) \log{(1-t^2)} - \log{\left( \prod \limits_{j=0}^{n-1}1-\mu_j t + (d-1)t^2 \right) }\\
& = -(|E|-|V|) \log{(1-t^2)} - \log{ \left(det(I-At+(d-1)t^2) \right)}
\end{align*}
Thus,
\begin{theorem}[Bass's determinant formula]
Let $G=(V,E)$ be a $d$-regular graph with adjacency matrix $A$, and let $N_k$ count the number of rooted, non-backtracking cycles of length $k$ in $G$. Then
$$\zeta_G(t) = \frac{ (1-t^2)^{|V|-|E|} } {det(I-At+(d-1)t^2)}$$
\end{theorem}

\end{document}